\pdfoutput=1

\documentclass{amsart}       % onecolumn (second format)

\RequirePackage{fix-cm}

\usepackage{graphicx}  % needed for figures
\usepackage{dcolumn}   % needed for some tables
\usepackage{bm}        % for math
\usepackage{amssymb}  
\usepackage{amsthm}
\usepackage{amsmath}
\usepackage{epstopdf}
\usepackage{tikz}
\usepackage{pgfplots}
\usepackage{mathtools}
\usepackage{amsaddr}
\usepackage{todonotes}

\let\mathds\mathbb

\usetikzlibrary{external}
\newcommand{\aspect}{\kappa}

\begin{document}

% The following information is for internal review, please remove them for submission

% the following line is for submission, including submission to the arXiv!!
%\hspace{5.2in} \mbox{Fermilab-Pub-04/xxx-E}

\title[Stability Results for Idealised Shear Flows]{Stability Results for Idealised Shear Flows on a Rectangular Periodic Domain}
\author{Holger Dullin and Joachim Worthington}
\address[A1,A2]{School of Mathematics and Statistics,
					Carslaw Building (F07),\\
			The University of Sydney, NSW 2006\\
	Email: Joachim.Worthington@sydney.edu.au (corresponding author), Holger.Dullin@sydney.edu.au}

\date{\today}

\begin{abstract}
We present a new linearly stable solution of the Euler fluid flow on a torus.
On a two-dimensional rectangular periodic domain $[0,2\pi)\times[0,2\pi / \aspect)$ for $\aspect\in\mathbb{R}^+$, the Euler equations admit a family of stationary solutions given by the vorticity profiles
$\Omega^*(\mathbf{x})= \Gamma \cos(p_1x_1+ \aspect p_2x_2)$. 
We show linear stability for such flows when $p_2=0$ and $\aspect \geq |p_1|$ (equivalently $p_1=0$ and $\aspect{|p_2|}\leq{1}$).
The classical result due to Arnold is that for $p_1 = 1, p_2 = 0$ and $\aspect \ge 1$ the stationary flow is {nonlinearly} stable via the energy-Casimir method.
We show that for $\aspect \ge |p_1| \ge 2, p_2 = 0$ the flow is linearly stable, but one cannot expect a similar nonlinear stability result. 
Finally we prove nonlinear instability for all equilibria satisfying $p_1^2+\aspect^2{p_2^2}>\frac{{3(\aspect^2+1)}}{4(7-4\sqrt{3})}$. 
The modification and application of a structure-preserving Hamiltonian truncation is discussed for the $\aspect\neq 1$ case. 
This leads to an explicit Lie-Poisson integrator for the truncated system.
\end{abstract}

\maketitle

\section{Introduction}
The study of nonlinear stability of stationary solutions of the 2D Euler equations was initiated by Arnold's seminal work  \cite{Arnold66} on the Lie-Poisson structure of the Euler equations. For a periodic domain Arnold considered a stationary solution with vorticity $\cos(x_1)$. This solution is nonlinearly stable in the ``energy-Casimir'' sense for domain size $\aspect\geq 1$ \cite{Arnold98}.  A detailed exposition of Arnold's method can be found in \cite{holm85}.
We will show that the stationary solution $\cos(mx_1)$ is linearly stable when $\aspect\geq m$ for any $m\in\mathbb{N}$. 
It may seem as though a scaling argument could give such a result; 
increasing the number of oscillations in the $x_1$-direction while decreasing the length of the domain in the $x_2$ direction.
However, we will show that nonlinear stability is lost in this process, allowing for instability on exponentially long timescales.
In the following we also consider the more general family of stationary solutions with vorticity $\Omega^*=\Gamma\cos(p_1x_1+\aspect p_2x_2)$ for $\mathbf{p}=(p_1,p_2)\in\mathbb{Z}^2$, $\Gamma\in\mathbb{R}$. 
The (in)stability of this family has been widely studied, e.g.\ \cite{friedlander97,Li00,butta10,meshalkin61,Beck12}. In particular \cite{Friedlander99} shows that for $p_2=0$, $p_1<\aspect$ there is linear instability. Classical methods like Rayleigh's inflection point criterion (see, e.g., \cite{drazin04,schmid12}) cannot give us stability in the present setting because
a shear flow on the torus is periodic, and a smooth periodic function always has inflection points.
We extend a number of existing results for $\aspect=1$ to $\aspect\neq 1$, and discuss the changes made by deforming the domain.

The Euler equations can be considered the inviscid limit of the Navier-Stokes equations under certain conditions, see, e.g.\ \cite{constantin95}.  Stability of stationary solutions to the Euler equations is used to predict and analyse the persistence of vortex structures in the Navier-Stokes equations. The idea is that unstable manifolds of steady solutions in the 2D Euler equations persist in the 2D Navier-Stokes problem for high value of the Reynolds number \cite{Li00,cichowlas05}. Such vortex structures are studied in the search for possible finite-time blowup solutions in the Navier-Stokes problem.

Written in terms of the vorticity $\Omega(\mathbf{x},t)$, the 2D incompressible Euler equations are 
\begin{equation}
\label{eq:originalsystem}
\frac{\partial \Omega}{\partial t} +u_1 \frac{\partial \Omega}{\partial x_1}+u_2\frac{\partial \Omega}{\partial x_2}=0,\;\;\;\;\; \frac{\partial u_1}{\partial x_1}+\frac{\partial u_2}{\partial x_2}=0.
\end{equation}
Here $\mathbf{x}=(x_1,x_2)^T$,  and $u_1 = \partial \psi/\partial x_2$, $u_2 =-\partial \psi/\partial u_1 $ are the velocity components in the $x_1$ and $x_2$ directions respectively,
and $\psi$ is the stream function related to the vorticity by $\Omega = \Delta \psi$.
We impose periodic boundary conditions $\Omega(\pi,x_2,t)=\Omega(-\pi,x_2,t)$ and $\Omega(x_1,{\pi}/{\aspect},t)=\Omega(x_1,-{\pi}/{\aspect},t)$ for $\aspect\in\mathbb{R}^+$.  Then there is a family of stationary solutions given by 
\begin{equation}
\label{eq:myequilibria}
	\Omega^*=2{\Gamma}\cos(p_1x_1+\aspect p_2x_2)
\end{equation}
for constant real $\Gamma$.
\begin{figure}
\includegraphics[width=0.68\textwidth]{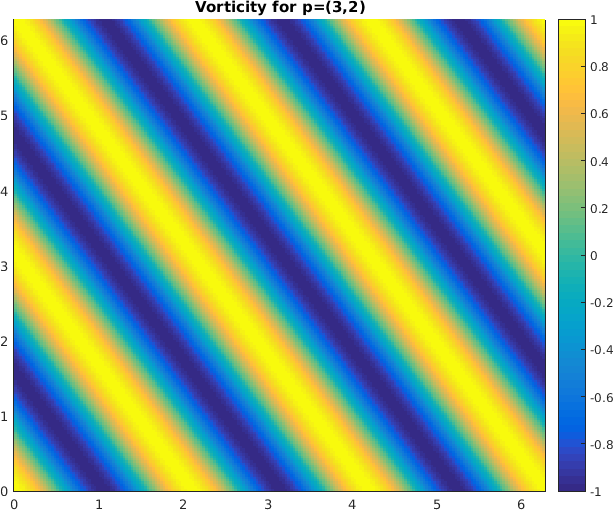}
\caption{The vorticity defined by \eqref{eq:myequilibria}. For this figure, $\mathbf{p}=(3,2)$, $\aspect=1$. The boundaries are periodic.}
\label{fig:velocity}
\end{figure}
The vorticity for this steady state is illustrated in Figure \ref{fig:velocity}.

In this paper we  study the linearisation of \eqref{eq:originalsystem} around the steady states \eqref{eq:myequilibria}. 
In Fourier space this leads to a linear operator that decomposes into subsystems or ``classes'' per \cite{Li00}.  
We present a new derivation of this class decomposition which uses the Poisson structure of the problem and leads naturally 
 to quadratic forms relevant for proving nonlinear stability. 
Our main stability result is based on the idea that if a region called the ``unstable ellipse'' in Fourier space only contains trivial lattice points $\mathbf k \in \mathbb Z^2$, 
then the equilibrium is linearly stable. A trivial lattice point satisfies $\mathbf k \parallel \mathbf p$ and does not contribute non-zero eigenvalues.
In addition instability is proved for most equilibria of the family \eqref{eq:myequilibria} by extending results from \cite{DMW16} to the case $\aspect \not = 1$.

We also modify the sine-bracket truncation of the Euler equation as described in \cite{Zeitlin90} for the $\aspect\neq 1$ case. 
This leads to an explicit Lie-Poisson integrator very similar to \cite{mclachlan93}. We illustrate our analytical results with some numerical experiments
using this structure preserving integrator.

\section{The Linearised 2D Periodic Euler Equations}
As the domain is periodic, we can expand the vorticity as $$\Omega(\mathbf{x},t)=\sum_\mathbf{k}\omega_\mathbf{k}(t)e^{i(k_1x_1+lk_2x_2)}$$ for $\omega_{\mathbf k}(t)\in\mathbb{C}$
where the sum is over all $\mathbf k \in \mathbb Z^2 \setminus \{ \mathbf 0 \}$. 
The condition $\omega_{\mathbf k}=\overline{\omega_{-\mathbf k}}$ is imposed for real-valued vorticity.
The dynamics of the Fourier coefficients $\omega_{\mathbf k}$ are then governed by a Hamiltonian system with the non-canonical ideal fluid Poisson bracket 
\begin{equation}
	\{f,g\}=\sum_{\mathbf{l},\mathbf{j}} \frac{\partial f}{\partial \omega_{\mathbf{l}}}\frac{\partial g}
	{\partial \omega_{\mathbf{j}}} {\aspect} (\mathbf{l}	\times \mathbf{j})\omega_{\mathbf{l}+\mathbf{j}}
	\label{eq:bracket}
\end{equation}
(where $\mathbf{a}\times\mathbf{b}=a_1b_2-a_2b_1$) 
and the Hamiltonian  is the kinetic energy of the fluid per unit area $\int |\nabla \psi|^2 dx_1 dx_2$ and hence
\begin{equation}
H=  \frac{1}{2}\sum_{\mathbf{k}\in \mathds{Z}^2\backslash \{\mathbf{0}\}} 
\frac{\omega_{+\mathbf{k}}\omega_{-\mathbf{k}}}{|\mathbf{k}|_{\aspect}^2} \,.
	\label{eq:untruncatedhamiltonian}
\end{equation}
Here $|\mathbf{k}|_{\aspect}=\sqrt{k_1^2+\aspect^2k_2^2}$ is the weighted norm for the rectangular domain.
For $\kappa=1$ \eqref{eq:bracket} matches the well-known ideal fluid Poisson bracket.
The (infinitely many) differential equations for the modes are
\begin{equation}
	\dot{\omega}_{\mathbf{k}}(t)= \{ \omega_{\mathbf k},  H \} =  
	\sum_{\mathbf{j}\in\mathds{Z}^2 \backslash \{\mathbf{0}\}}\frac{\aspect\mathbf{k} \times \mathbf{j}}
	{\;|\mathbf{j}|_{\aspect}^2}\omega_{-\mathbf{j}}\omega_{\mathbf{k}+\mathbf{j}}.
\label{eq:fourierodes}
\end{equation}

We now linearise about the steady state \eqref{eq:myequilibria}, given in Fourier space by 
\begin{equation}
\omega_{\mathbf k}^* =
\begin{cases} 
   \Gamma & \text{ if } \mathbf{k}=\pm \mathbf{p}, \\
   0       & \text{otherwise.}
  \end{cases}
\end{equation}

\begin{figure}
\centering
\includegraphics{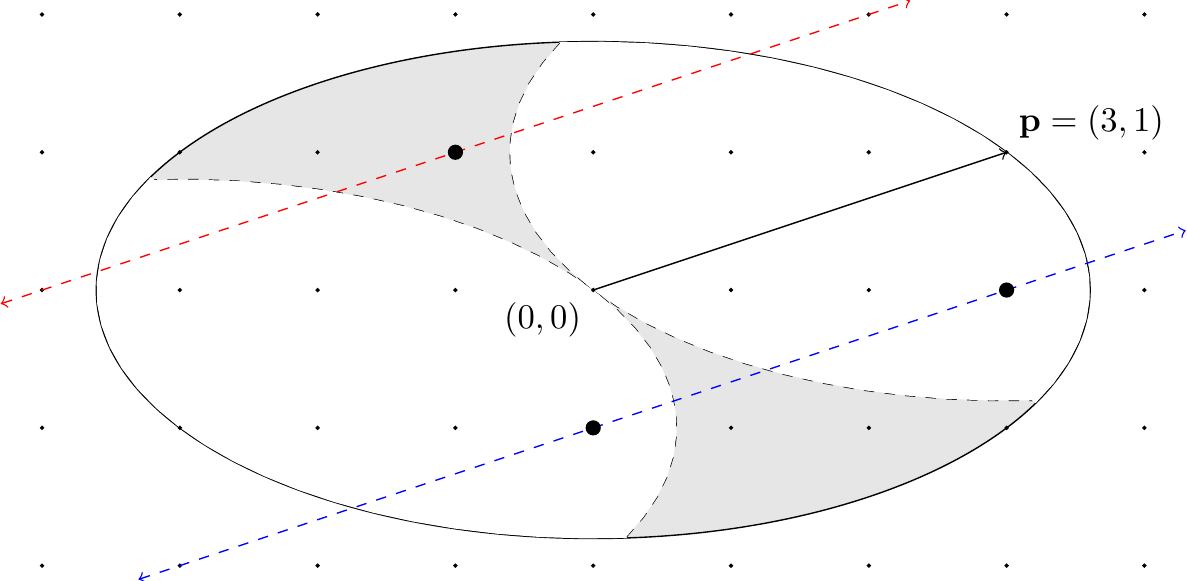}

\caption{The class structure of the problem. Each lattice point $\mathbf{k}$ corresponds to a Fourier mode $\omega_\mathbf{k}$, the dynamics of which depend only on $\omega_{\mathbf{k}+\mathbf{p}}$ and $\omega_{\mathbf{k}-\mathbf{p}}$ in the linearised system. 
The modes that lie along a line parallel to $\mathbf{p}$ belong to the same ``class''. Some example classes are shown here by dashed straight lines. The indicated ``unstable ellipse'' (cf {unstable disc} for $\aspect=1$ \cite{Li00}) 
with centre $\mathbf{0}$ and axes $|\mathbf{p}|_{\aspect}$ and $|\mathbf{p}|_{\aspect}/\aspect$ delineates the classes. Classes that do not intersect the ellipse are linearly stable. Classes that intersect the disc are linearly unstable. In this figure $\mathbf{p}=(3,1)$, $\aspect=2$. }
\label{fig:classes}
\end{figure}

The following calculation can be thought of as a Fourier space version of the first few steps of the stability ``algorithm'' described in \cite{holm85}.
When linearising a system with a noncanonical Poisson bracket the gradient of the Hamiltonian does not necessarily vanish at an equilibrium. 
This is due to the presence of Casimir functions, which have vanishing Poisson bracket with any function. The bracket \eqref{eq:bracket} has a
family of Casimirs $C_n$ given by 
\begin{equation}
     C_n=\sum_{\sum \mathbf k_i=\mathbf 0} \omega_{\mathbf k_1}\omega_{\mathbf k_2}...\omega_{\mathbf k_n}, \;\;n\in\mathbb{N} 
\label{eq:casimirs}
\end{equation}
obtained from the conserved quantities $\int \Omega^n dx_1 dx_2$ \cite{Arnold98,Zeitlin90}.
 A particular linear combination of gradients of $H$ and gradients of $C_n$ vanishes at an equilibrium. In our case define
 \begin{equation}
     F=H-\frac{1 }{2|\mathbf{p}|_{\aspect}^2}C_2 =  
       - \frac{1}{2}\sum_{\mathbf{k \neq 0}}  \rho_{\mathbf k}     |\omega_\mathbf{k}|^2, \quad
      \rho_{\mathbf k} = \frac{1}{|\mathbf{p}|_{\aspect}^2} - \frac{1}{|\mathbf{k}|_{\aspect}^2},
\end{equation}
such that $F_{\Omega^*}=0$, $\nabla F | _{\Omega^*}=\mathbf{0}$. 
Denote by $P$ the matrix of the Lie-Poisson structure with entries that are linear in $\omega$. 
The vector field can be written as $P \nabla H$, and for the Casimirs we have $P \nabla C_n = 0$. 
Assuming that the $\nabla C_n$ span the kernel of $P$ we can write $\nabla H|_{\Omega^*}$ as a linear combination 
of $\nabla C_n|_{\Omega_*}$, and the difference leads to $F$ above. Now linearising 
$P \nabla H $ and $P \nabla C_n$ and combining it with $\nabla F = \mathbf{0}$ gives
the linearised matrix $P^* D^2 F|_{\Omega^*}$ where $P^* = P|_{\Omega_*}$.
Thus we have shown that the linearisation at $\Omega^*$ can be found by evaluating the Poisson bracket with $F$ at $\Omega^*$. 
Hence (again writing $\omega_\mathbf k$ for the linearisation about $\Omega^*$ from now on) 
\begin{equation}
   \dot \omega_\mathbf k \!= \! \left. \{ \omega_\mathbf k, F \}\right|_{\Omega^*} \!\! =
     \sum_{\mathbf{k \not = 0}} \rho_{\mathbf k} ( \mathbf{ k \times j}) \omega_{-\mathbf j} \Gamma (\delta_{\mathbf{k+j,p}} + \delta_{\mathbf{k+j,-p}} ).
\end{equation}
This shows that the linearised equations decouple into subsystems or ``classes'' \cite{Li00}: 
the mode $\mathbf k$ only couples to the modes $\mathbf{ k \pm p}$, and hence a subsystem 
consists of modes $\mathbf{ a} + k \mathbf p$ for $k \in \mathbb Z$ and any fixed $\mathbf a \in \mathbb Z^2$,
illustrated in Fig.~\ref{fig:classes}.

To make the splitting into subsystems explicit write 
\begin{equation} \label{eqn:Fa}
   F = -\sum_{\mathbf a \in A} F_{\mathbf a}, \quad
    F_{\mathbf a} =  \frac{1}{2}\sum_{k \in \mathbb Z} \rho_k \omega_k^2
\end{equation}
where 
$A = \{\mathbf a \in \mathbb Z^2 \; | \; -|\mathbf p|_\aspect^2 < 2 ( a_1 p_1 + \aspect^2 a_2 p_2) \le |\mathbf p|_\aspect^2 \} $ and
\begin{equation} \label{eq:rhodef}
   \rho_k = \rho_{\mathbf a + k \mathbf p}, \quad
   \omega_k = \omega_{\mathbf a + k \mathbf p} .
\end{equation}

Then for each $\mathbf{a}\in A$ there is a linearised subsystems which is Hamiltonian with respect to the Poisson bracket
\begin{equation}
\label{eq:linbracket}
\{f,g\}_{\mathbf a}=\alpha \sum_{k\in\mathbb{Z}} \left [ \frac{\partial f}{\partial \omega_{k}}\frac{\partial g}{\partial \omega_{k+1}}-\frac{\partial f}{\partial \omega_{k+1}}\frac{\partial g}{\partial \omega_{k}}\right ]
\end{equation}
with Hamiltonian $F_\mathbf{a}$, where 
\begin{equation}
	\alpha = \alpha(\mathbf{a},\mathbf{p})=\Gamma  \aspect \mathbf{a}\times\mathbf{p}.
\end{equation}
It is important to note that if $\mathbf{a}$ and $\mathbf{p}$ are parallel, 
$\alpha(\mathbf{a},\mathbf{p})=0$.
Such ``trivial'' classes cannot contribute instability as the associated ODEs are $\dot{\omega}_k=0$ for all $k$.

There are some symmetries in the linearised problem. Define the triple $(\mathbf{a},\mathbf{p},\aspect)$ to refer to the class led by $\mathbf{a}$ with parameters $\mathbf{p}$, $\aspect$. Then $(\mathbf{a},\mathbf{p},\aspect)$,  and $(\pm n\mathbf{a},n\mathbf{p},\aspect)$ have the same spectrum for any nonzero $n\in\mathbb{Z}$. This is clear from \eqref{eq:rhodef} and \eqref{eq:MainMatrixG}. Thus if the stationary solution with parameters $\mathbf{p}$, $\aspect$ is unstable, so is the stationary solution with $ n\mathbf{p}$, $\aspect$ for nonzero $n\in\mathbb{Z}$. The same conclusion cannot be made for stable stationary classes, as there will be classes in the $n\mathbf{p}$, $\aspect$ problem that do not have a corresponding class in the $\mathbf{p}$, $\aspect$ problem. 

Similarly the $((a_1,a_2),(p_1,p_2),\aspect)$ and $((a_2,a_1),(p_2,p_1),\aspect^{-1})$ problems will have the same associated spectrum. For every class in the $(p_1,p_2),\aspect$ problem there is a corresponding class in the $(p_2,p_1),\aspect^{-1}$ problem with identical spectrum, so if one is stable (respectively unstable), the other is also stable (respectively unstable).

Define the {\em unstable ellipse} as
$D_\mathbf{p}=\{\mathbf{x}\in \mathbb R^2 \; | \; |\mathbf{x}|_{\aspect}<|\mathbf{p}|_{\aspect}\;\}$, as shown in Fig.~\ref{fig:classes}. We associate lattice points in this region with the corresponding Fourier modes.
Note that $D_\mathbf{p}$ does \emph{not} include the boundary.
Note that $\rho_\mathbf{k}<0$ if and only if $\mathbf{k}\in D_\mathbf{p}$.
Then we can say that the coefficients $\rho_{\mathbf k}$ 
are all non-negative if and only if $\mathbf{a}\notin D_\mathbf{p}$.

Explicitly the ODE for the  bi-infinite vector of Fourier coefficients $\omega=(...,\omega_{-1},\omega_0,\omega_1,\omega_2,...)$ in the subsystem $F_\mathbf{a}$ is 
\begin{equation}
\dot{\mathbf{\omega}}=\alpha M \mathbf{\omega}
 \label{eq:linode}
\end{equation}
 where $\alpha$ is constant and
 {\footnotesize
\begin{equation}
\label{eq:MainMatrixG}
M = \begin{pmatrix}
	\ddots & \vdots & \vdots & \vdots & \vdots & \vdots & \ddots \\	
	\cdots & 0 &	\rho_{-1} & 0 & 0  &0  & \cdots  \\
	\cdots & -\rho_{-2} &	0 & +\rho_{0} & 0  & 0  &\cdots  \\
	\cdots & 0 &	-\rho_{-1 } & 0 & +\rho_{1} &  0  &\cdots  \\
	\cdots & 0 &	0 & -\rho_{0} & 0 & \rho_{2} & \cdots \\
	\cdots & 0 &	0 & 0 & -\rho_{1} &  0  &\cdots  \\
	\ddots & \vdots &	\vdots & \vdots & \vdots  & \vdots &  \ddots 
	\end{pmatrix}.
\end{equation} }

\section{Stability Analysis - Stable Values of $\mathbf{p}$}

\label{sec:stable}

The question of stability is now that of calculating the spectrum of $M$. Note that $M$ is an Hamiltonian operator, so $\lambda\in\sigma(M)$ implies $-\lambda,\bar{\lambda},,-\bar{\lambda}\in\sigma(M)$. 
The operator $M$ has continuous spectrum on some part of the imaginary axis, and may have discrete eigenvalues off the imaginary axis.

If $\mathbf a \notin  { D_\mathbf p}$, then $\rho_k \geq 0$ for all $k$. Define 
\begin{equation}
	J=\begin{pmatrix}
	\ddots & \vdots & \vdots & \vdots  & \vdots & \ddots \\	
	\cdots & 0 & +1 & 0 & 0    & \cdots  \\
	\cdots & -1  &	0 & +1  & 0  &\cdots  \\
	\cdots & 0 &	-1 & 0 & +1 &\cdots  \\
	\cdots & 0 &	0 & -1 & 0 &\cdots \\
	\ddots & \vdots &	\vdots & \vdots   & \vdots &  \ddots 
	\end{pmatrix}.
\end{equation}
and $T=\text{diag}( ...,\sqrt{\rho_{-1}},\sqrt{\rho_0},\sqrt{\rho_1},...)$. $J$ is a real skew-symmetric operator and $T$ is a real symmetric operator. Then $M=JT^2$. Define $\tilde{M}=TJT$. Then $i\tilde{M}$ is Hermitian and so $\tilde{M}$ has only imaginary spectrum. As $J,T$ are bounded operators, $M$ and $\tilde{M}$ have the same spectrum and so $M$ has no imaginary eigenvalues and cannot contribute instability. Thus all such classes are linearly stable. Geometrically, such classes have all lattice points lying {outside} the unstable ellipse.   
This statement generalises \cite{Li00,shvidkoy03} to the case $\aspect \not = 1$. 

\begin{figure}
\centering
\includegraphics{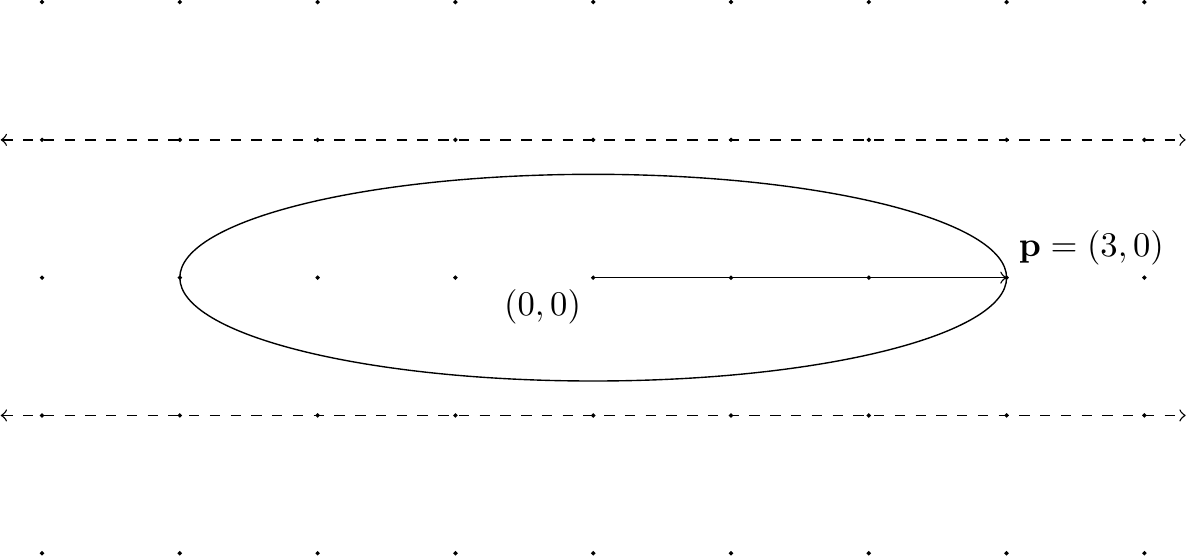}

\caption{The unstable ellipse for a set of parameters leading to a linearly stable stationary solution. Here $\mathbf{p}=(3,0)$ and $\aspect=4$. If $\mathbf{p}=(p_1,0)$ and $\aspect\geq p_1$, there are no lattice points inside the unstable ellipse that are not parallel to the vector $\mathbf{p}$. Thus no classes contribute linear instability and thus there is linear stability.}
\label{fig:stableclasses}
\end{figure}

We can further characterise the stable spectrum. Note that  $\lim_{k\to \pm \infty} \rho_k = 1/|\mathbf p|_\aspect^2 = const$, and hence the Hermitian operator is unitarily equivalent 
to a Jacobi operator that can be written as $J_0 + P$ where $J_0=1/ |\mathbf p|_\aspect^2 J$ and the perturbation 
$P$ is a compact operator whose coefficients decay as $1/k^2$ for $k \to \pm \infty$. By Weyl's theorem on the spectrum of compact perturbations of operators, see, e.g. \cite{reed78},
the essential spectrum of $J_0 + P$ is equal to the essential spectrum of $J_0$. The constant coefficient recursion relation $J_0 v = \lambda v$ 
has solution $v_k = \mu^k$ where $\mu$ are the roots of the characteristic equation $1 + \mu^2 = ( \lambda  |\mathbf p|_\aspect^2/ 2) \mu$,
and hence bounded solutions exists for $\mu \in [ - 2 , 2]/ |\mathbf p|_\aspect^2$. %\todo{is bounded enough?}
For more details on this argument in the case $\aspect =1$ see \cite{Li04}.

This leads to our main result: if $\mathbf{p}=(p_1,0)$ and $\aspect \geq |p_1|$, the stationary solution 
\begin{equation}
\label{eq:stableequilibria}
	\Omega^*=2{\Gamma}\cos(p_1x_1)
\end{equation}
is linearly stable. This is because $\rho_k<0$ if and only if $|\mathbf{a}+k\mathbf{p}|_{\aspect}<|\mathbf{p}|_{\aspect}=|p_1|$. If $\aspect\geq |p_1|$, this condition implies $a_2=0$ so so $\mathbf{a}=(a_1,0)$. But then $\mathbf{a}$ and $\mathbf{p}$ are parallel and thus $\alpha(\mathbf{a},\mathbf{p})=0$. Thus all classes either do not intersect the unstable ellipse, or they satisfy $\alpha(\mathbf{a},\mathbf{p})=0$ and so there are no classes that contribute instability. An example of this is illustrated in Figure \ref{fig:stableclasses}. Analogously, if $\mathbf{p}=(0,p_2)$ and $\aspect\leq \frac{1}{|p_2|}$ the corresponding stationary solution is linearly stable.  

We emphasise that this linear stability is optimal and cannot be extended to nonlinear energy-Casimir stability in the sense of \cite{Arnold98}.  The general procedure for showing nonlinear stability in the energy-Casimir sense is discussed in detail in \cite{holm85} (see also \cite{morrison98}).

To understand the nonlinear instability we restrict the quadratic form $F$ to the complement of the gradients of the Casimirs. This allows us 
to consider stability while locally fixing the values of the Casimirs. 

The linear approximation of the Casimirs $C_n$ at the equilibrium $\Omega^*$ 
is given by a linear combination of the modes $\omega_{k \mathbf p}$ where $k = -n, -n+2, \dots, n-2, n$. 
Thus the span of the linear approximations of all Casimirs $C_1,C_2,...$ is \emph{all} linear combinations of the modes $\omega_{k\mathbf{p}}$ for $k\in\mathbb{Z}$.

In the case $\aspect=1$ the Hessian $D^2 F$ has eigenvalues zero originating from $\rho_{\pm\mathbf p} = 0$. These neutral directions are removed by fixing the Casimirs,
and energy-Casimir stability follows.
This is the Fourier space version of Arnold's theorem (for the original, see section 4 of \cite{Arnold98}). 

For our new linearly stable cases with $\aspect>|p_1|>1$, the quadratic form $F$ is indefinite, specifically at lattice points $\omega_\mathbf k$ where $\mathbf k \in D_\mathbf{p}$. This implies they are of the form $\mathbf{k}=(k_1,0),\;|k_1|<|p_1|$ so $\rho_\mathbf k < 0$.
One cannot achieve definiteness of $F$ by restricting to fixed values of the Casimirs for these modes, as they are not at integer multiples of $\mathbf{p}$. That is, the second variation cannot be made definite by fixing the Casimirs, so we cannot find the convexity estimates required by \cite{holm85}. Our argument is equivalent to considering the most general 
linear combination of Casimirs to construct $F$ as suggested in \cite{holm85}, and then to observe that even this most general $F$ is not definite for the case at hand.

We thus cannot conclude nonlinear energy-Casimir stability for these linearly stable steady states.

\section{Stability Analysis - Unstable Values of $\mathbf{p}$}

In addition to the stability results, we can make some conclusions for unstable steady states. Take a Galerkin-style truncation of \eqref{eq:linode}. Select some $m,n\in\mathbb{Z}$ with $m<n$, and project the class to the modes $\omega_m,\omega_{m+1},...,\omega_{n-1},\omega_n$ by setting all other modes to zero. 
The ODEs for this truncated system are
\begin{equation}
\dot{\mathbf{\omega}}=\alpha A_m^n \mathbf{\omega}
\end{equation}
 where
 {\footnotesize
\begin{equation}
\label{eq:MainMatrixA}
A_m^n= \begin{pmatrix}
		0 & +\rho_{m+1} & 0  & \cdots &  0&  0 & 0 \\
		-\rho_{m} & 0 & +\rho_{m+2} &  \cdots  &  0& 0 &  0 \\
		0 & -\rho_{m+1} & 0 &  \cdots &  0 & 0 & 0 \\
		\vdots & \vdots & \vdots &  \ddots &  \vdots & \vdots & \vdots \\
		0 & 0 & 0 &  \cdots &  0 &  +\rho_{n-1} & 0  \\
		0 & 0 & 0 &  \cdots &  -\rho_{n-2} &  0 & +\rho_{n}  \\
		0 & 0 & 0 &  \cdots &  0 &  -\rho_{n-1} & 0
	\end{pmatrix}
\end{equation} } 
and $\omega=(\omega_{m},\omega_{m+1},...,\omega_{n-1},\omega_n)$. We make the following conjectures based on numerical observations:
\begin{enumerate}
	\item if $\rho_k<0$ for exactly one value $k\in\mathbb{Z}$, $A_m^n$ has a pair of real eigenvalues;
	\item if $\rho_k<0$ for two values of $k$, $A_m^n$ has four non-imaginary eigenvalues; either a complex quadruplet of nonimaginary eigenvalues or two pairs of real eigenvalues.
\end{enumerate}

Figure \ref{fig:classes} gives a geometric interpretation of these cases. Classes with $\mathbf{a}$ outside the unstable ellipse will be stable. Classes with $\mathbf{a}$ in the shaded region in Figure \ref{fig:classes} will lead to a pair of real eigenvalues. Classes with $\mathbf{a}$ in the unstable ellipse but outside the shaded area will lead to four non-imaginary eigenvalues. 

We can prove Case 1 with a proof similar to \cite{DMW16}. However, a concern is that the positive real eigenvalue may vanish in the limit $n\to\infty$, $m\to -\infty$. We therefore need to show that there is some positive finite lower bound on a real eigenvalue.

The result is as follows: if $\rho_0<0$, $\rho_k>0$ for all $k\neq 0$, and $\rho_0+\rho_2<0$, then there is some eigenvalue $\lambda$ of $A_m^n$  such that 
\begin{equation}
	\label{eq:lowerbound}
	\lambda\geq\lambda^*=\sqrt{-\rho_1(\rho_0+\rho_2)}>0.
\end{equation}
As this bound does not depend on $n$ and $m$, the positive real eigenvalue will persist in the limit $n\to\infty$, $m\to -\infty$.

To prove this, consider the characteristic polynomials $P_m^n(x)=\text{det}(xI-A_m^n)$ for truncation values $m,n$. By a well-known result  for the determinants of tridiagonal matrices, 
\begin{align} 
P_m^n(x)&=xP_m^{n-1}(x)+\rho_n\rho_{n-1}P_{m}^{n-2}(x) \label{eq:charrecursion}\\
 &=xP_{m+1}^{n}(x)+\rho_{m}\rho_{m+1}P_{m+2}^{n}(x).
\end{align}

Assume $\rho_0<0$, $\rho_k>0$ for all $k\neq 0$, and $\lambda^*\in\mathbb{R}$. It is straightforward to check that $P_{0}^{2}(\lambda^*)=0$, $P_{0}^{3}(\lambda^*),P_{-1}^{2}(\lambda^*),P_{-1}^{3}(\lambda^*)<0$. Then by induction on \eqref{eq:charrecursion}, $P_{m}^{n}(\lambda^*)\leq 0$ for all $m\leq 0$ and $n\geq 2$. The leading order term of $P_{m}^{n}(x)$ is $x^{m+n+1}>0$ so $\lim_{x\to\infty} P_{m}^{n}(x)>0$. Therefore by the intermediate value theorem, for all $m\leq 0,\;n\geq 2$ there is some $\lambda\geq \lambda^*$ such that $P_{m}^{n}(\lambda)=0$. Thus for any  truncation values an eigenvalue with a non-vanishing positive real part exists. This eigenvalue persists in the limit, the untruncated linearised problem.

\begin{figure}
\centering
\includegraphics{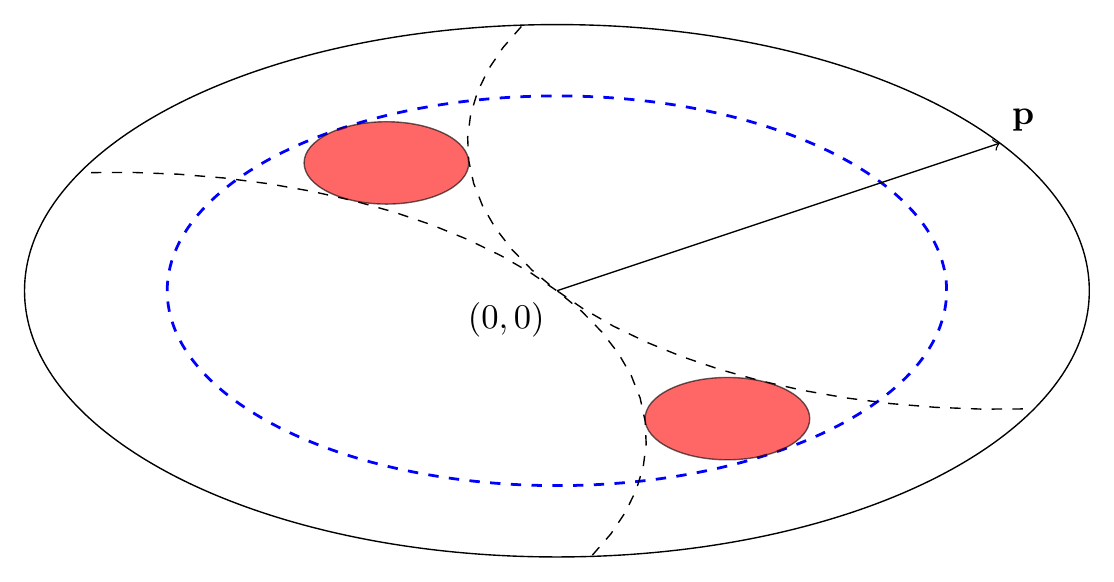}

\caption{The dashed blue ellipse has origin $(0,0)$ and axes $(\sqrt{3}-1)|\mathbf{p}|_{\aspect}$ and $(\sqrt{3}-1)\frac{ |\mathbf{p}|_{\aspect}}{\aspect}$. If $\mathbf{a}$ is in the interior of this ellipse, $\rho_0+\rho_2<0$ and $\rho_0+\rho_{-2}<0$. The red shaded regions lie wholly within the intersection of this condition and the condition $\rho_0<0$, $\rho_1,\rho_{-1}\geq 0$ as shown in Figure \ref{fig:classes}. The requirements of the instability proof are satisfied if there is a lattice point in this region, which is true if \eqref{eq:bound} is satisfied.}

\label{fig:Conditions2}
\end{figure}

We need to determine which values of $\mathbf{p}$ admit a value of $\mathbf{a}$ such that $\lambda^*$ is real and positive. A sufficient condition is the existence of an integer lattice point in the shaded region indicated in Figure \ref{fig:Conditions2}. Algebraically, this is possible if
\begin{equation}
	|\mathbf{p}|_{\aspect}>\frac{\sqrt{3(\aspect^2+1)}}{{2}(2-\sqrt{3})}.
	\label{eq:bound}
\end{equation}
Compare this with equation (4.5) in \cite{DMW16}. For a fixed finite $\aspect\in\mathbb{R}^+$, this condition is satisfied for all but finitely many values of $\mathbf{p}\in\mathbb{Z}^2$. 

For all $\mathbf{p}$ and $\aspect$ satisfying this condition this implies \emph{nonlinear instability} per the spectral gap theorem of \cite{friedlander97} and \cite{shvidkoy03}. The extrapolation from linear instability to nonlinear instability is discussed in \cite{DMW16}.

Note that this bound is not sharp. 
Numerical evidence suggests that in fact any $\mathbf{p}=(p_1,p_2)$ and $\aspect$ \emph{not} satisfying $p_2=0, \aspect>|p_1|$ or $p_1=0,\aspect |p_2|<1$ will lead to an unstable stationary solution.

 \section{Non-Imaginary Spectrum}
 
  \begin{figure}
 \includegraphics[width=0.7\textwidth]{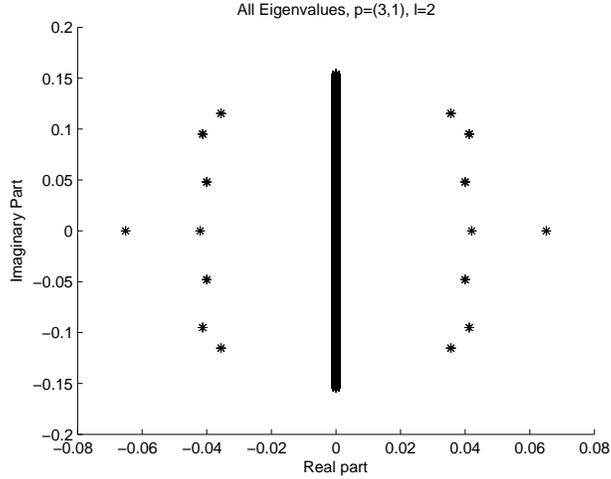}
 \caption{The full spectrum for $\mathbf{p}=(3,1)$, $l=2$, $\Gamma=1$. The truncation $m=-100$, $n=100$ was used for this figure.}
 \label{fig:eigenvalues2}
 \end{figure}

\begin{figure}
 \includegraphics[width=0.8\textwidth]{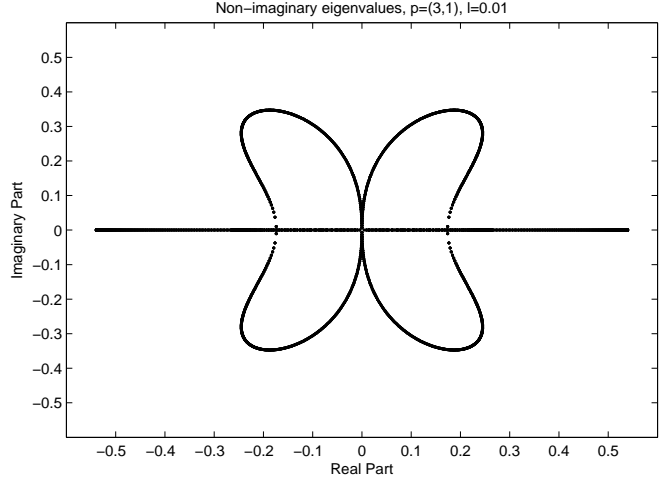}
 \caption{The non-imaginary spectrum for $\mathbf{p}=(3,1)$, $\aspect=0.01$, $\Gamma=1$. Truncation values used $-m=100=n$.}
 \label{fig:eigenvalues3}
 \end{figure}
 
 We describe a conjecture on the discrete spectrum. Write $\nu$ for the number of points $(x,y)\in\mathbb{Z}^2$ satisfying $x^2+\aspect^2y^2<|\mathbf{p}|_{\aspect}^2$, ie the number of lattice points in the unstable ellipse. From our hypotheses about the number of eigenvalues in an unstable class, we conclude that for $\mathbf{p}=(p_1,p_2)$, the number of discrete  eigenvalues counting multiplicity is 
 \begin{equation}
|\sigma_{\text{ess}}|=2(\nu-2\text{gcd}(p_1,p_1)+1).
 \end{equation}
 This is based on numerics for a range of parameter values. For $\aspect=1$ and $p_1,p_2$ coprime, this matches the upper bound on the number of discrete eigenvalues proved in \cite{Li04}.  For instance, in Figure \ref{fig:classes} we can count $16$ interior points. In Figure \ref{fig:eigenvalues2}, one can observe 16 discrete eigenvalues, each of which has multiplicity 2.

 In the limits $\aspect\to 0 / \infty$ and $p_1,p_2\neq 0$, $\nu\to\infty$ and therefore $|\sigma_{\text{ess}}|\to \infty$. In these cases the discrete spectrum appears to converge to some continuous spectrum. Figure  \ref{fig:eigenvalues3} shows the non-imaginary spectrum for $\mathbf{p}=(3,1)$ and a small value of $\aspect$.

\section{Structure Preserving Truncation and a Lie-Poisson Integrator}

For the Euler equation in a square domain $\aspect=1$,  Zeitlin \cite{Zeitlin90,Zeitlin05} and  Fairlie and Zachos \cite{fairlie89} describe a structure preserving Hamiltonian truncation of the Poisson bracket \eqref{eq:bracket}. The idea is to consider only the Fourier coefficients $\omega_\mathbf{a}$ where $\mathbf{a}\in\mathcal{D}_N=\{(a_1,a_2)| -N\leq a_1,a_2\leq N\}$ for some truncation value $N\in\mathbb{N}$. All other Fourier modes are set to zero. The ``sine bracket'' on these modes then is
\begin{equation}
	\{f,g\}=\sum_{\mathbf{k},\mathbf{l}} \frac{\sin{\varepsilon \mathbf{k}\times\mathbf{l}}}{\varepsilon} 
	\frac{\partial f}{\partial \omega_\mathbf{k}}\frac{\partial g}{\partial \omega_\mathbf{l}} \omega_{\mathbf{k}+\mathbf{l} \text{ mod } 2N+1}.
	\label{eq:zeitlinbracket}
\end{equation}
Here  $\varepsilon=\frac{2\pi}{2N+1}$, and the mode number ${k+l}$ is taken modulo $2N+1$ so $\mathbf{k}+\mathbf{l}\in\mathcal{D}_N$.  
As $N\to \infty$, this bracket approaches \eqref{eq:bracket} with $\aspect=1$ in a natural way.

This structure-preserving truncation conserves $2N+1$ Casimirs corresponding to those defined for the full system by \eqref{eq:casimirs} (see \cite{Zeitlin05} for an explicit construction of these). Considered with the Hamiltonian 
\begin{equation}
H=  \frac{1}{2}\sum_{\mathbf{k}\in \mathcal{D}_N\backslash \{\mathbf{0}\}} \frac{\omega_{+\mathbf{k}}\omega_{-\mathbf{k}}}{|\mathbf{k}|_{1}^2} 
	\label{eq:zeitlinhamiltonian}
\end{equation}
this is a finite dimensional approximation to the nonlinear Euler equations with some useful mathematical properties. Further discussion of this truncation is given in \cite{engo01,scovel1996survey}, and the theoretical background is discussed in \cite{hoppe89,hoppe91}.

For general $\aspect\in \mathbb{R}^+$, the only required change is  to introduce the factor  $\aspect$ in the Poisson bracket
\begin{equation}
	\{f,g\}=\sum_{\mathbf{k},\mathbf{l}} \aspect\frac{\sin{\varepsilon \mathbf{k}\times\mathbf{l}}}{\varepsilon} 
	\frac{\partial f}{\partial \omega_\mathbf{k}}\frac{\partial g}{\partial \omega_\mathbf{l}} \omega_{\mathbf{k}+\mathbf{l} \text{ mod } 2N+1},
	\label{eq:zeitlinbracket2}
\end{equation}
and use the $\kappa$ weighted norm in the Hamiltonian 
\begin{equation}
H=  \frac{1}{2}\sum_{\mathbf{k}\in \mathcal{D}_N\backslash \{\mathbf{0}\}} 
\frac{\omega_{+\mathbf{k}}\omega_{-\mathbf{k}}}{|\mathbf{k}|_{\aspect}^2}. 
	\label{eq:zeitlinhamiltonian2}
\end{equation}
Then \eqref{eq:zeitlinbracket2} approaches \eqref{eq:bracket} for $N\to\infty$ for all values of $\aspect$.  The new bracket \eqref{eq:zeitlinbracket2} conserves the same $2N+1$ Casimirs as \eqref{eq:zeitlinbracket}. 

A practical application of such a bracket is the development of a Lie-Poisson integrator which conserves all Casimirs. A Lie-Poisson integrator for the sine-bracket truncated Euler equations was developed in \cite{mclachlan93}. As the new bracket \eqref{eq:zeitlinbracket2} is identical up to a constant factor, the same integrator can be used for general $\aspect$ with only minor alterations. Such an integrator was used to numerically calculate the dynamics of $\Omega$ for various values of $\mathbf{p}$ and $\kappa$ to corroborate the results of this paper.

Of interest is that the same energy-Casimir argument as in Section \ref{sec:stable} can be followed in the sine-bracket truncated system, as there is an analogous set of Casimirs. However, the ``wrapping'' operation leads to a surprising result. If $\mathbf{p}=(p_1,0)$, $\kappa\geq p_1$, and $2N+1$ and $p_1$ are coprime, then the sine-bracket truncated system is nonlinearly stable in the energy-Casimir sense. This occurs as the wrapping of the bracket around the modes in $\mathcal{D}_N$ means that the set of Casimirs $C_n$ will now fix the values of all modes $\omega_{(a,0)}$ for all $-N<a<N$. Such a result is artificial in the sense that it only works for special truncation values $N$, but does give insight into the use of the sine-bracket truncation. Care must be taken when selecting truncation values for the finite-dimensional truncation.

\section{Conclusion}

We have shown that for sufficiently narrow domains, 
shear flows of the form $\cos(p_1x_1)$ can be linearly stable on a periodic domain. 
This result has also been shown to be optimal in the sense that it cannot be extended to energy-Casimir stability. For most other values of $\mathbf{p}$ we have demonstrated nonlinear instability of the steady flow $\cos(p_1x_1+\aspect p_2x_2)$.

We have also outlined an extension of the well-known sine-bracket truncation for the Euler equations to a non-square periodic domain. This leads to a straightforward development of a Lie-Poisson integrator. This is particularly useful for finite-dimensional numerical analysis.

Further work is necessary to clarify the connection between these results and existing work, especially in the context of the Euler problem as the inviscid limit of the Navier-Stokes problem. Another interesting avenue of future research could be a careful analysis of the limit $\aspect\to \infty / 0$, and any possible physical interpretations of this model. The bound \eqref{eq:bound} may also be sharpened.

It is possible that  approaches similar to those used here could shed some light on the study of stationary solutions of the Euler equations on the sphere. This has obvious geophysical applications, especially for the rotating sphere. A structure-preserving  truncation analogous to the one discussed in this paper exists for the rotating sphere, in both the Euler and Navier-Stokes problems \cite{Zeitlin04}; this may prove useful in further study.

\bibliographystyle{amsplain}
\bibliography{prf_preprint}

\end{document}